\documentclass[12pt]{article}
\usepackage{graphicx}
\usepackage{amssymb, cite}
\usepackage{setspace} 
\usepackage[hmargin=1.15in,vmargin=1.15in]{geometry} 

\newcommand{\be}{\begin{equation}}
\newcommand{\ee}{\end{equation}}
\newcommand{\bqn}{\begin{eqnarray}}

\newcommand{\eqn}{\end{eqnarray}}

\newcommand{\bd}{\begin{description}}
\newcommand{\ed}{\end{description}}

\newtheorem{stat}{}[section]

\def\bs{\begin{stat}}
\def\es{\end{stat}}
\def\ben{\begin{enumerate}}
\def\een{\end{enumerate}}

\def\bp{\noindent{\bf Proof}  \ \ \ }
\newcommand{\ep}{\hfill $\square$}

\makeatletter
\@addtoreset{equation}{section}

\makeatother

\begin{document}

\begin{center}

{\large {\bf PACKING  3-VERTEX PATHS
 IN 2-CONNECTED GRAPHS}}
\\[4ex]
{\large {\bf Alexander Kelmans}}
\\[2ex]
{\bf University of Puerto Rico, San Juan, Puerto Rico}
\\[0.5ex]
{\bf Rutgers University, New Brunswick, New Jersey}
\end{center}

\begin{abstract}
Let $v(G)$ and $\lambda (G)$ be the number 
of vertices and the maximum number of disjoint
3-vertex paths in $G$, respectively.
We give a construction that provides 
infinitely many   2-connected, cubic, bipartite, and planar  graphs such that  $\lambda (G) <   \lfloor v(G)/3 \rfloor $.
\\[1ex]
\indent
{\bf Keywords}: cubic, bipartite, planar,  
$\Lambda $-packing,  $\Lambda $-factor.
 
\end{abstract}

\section{Introduction}

\indent

We consider undirected graphs with no loops and 
no parallel edges. All notions and facts on graphs, that are  
used but not described here, can be found in \cite{BM,D,Wst}.
\\[1ex]
\indent
Given graphs $G$ and $H$, 
an $H$-{\em packing} of $G$ is a subgraph of $G$ 
whose components are isomorphic to $H$.
An $H$-{\em packing} $P$ of $G$ is called 
an $H$-{\em factor} if $V(P) = V(G)$. 
The $H$-{\em packing problem}, i.e. the problem of 
finding in $G$ an $H$-packing, having the maximum 
number of vertices, turns out to be $NP$-hard if $H$ is 
a connected graph with at least three vertices \cite{HK}.
Let $\Lambda $ denote a 3-vertex path.
In particular, the $\Lambda $-packing problem
is $NP$-hard. Moreover, this problem  remains 
$NP$-hard even for cubic graphs \cite{K1}.

Although the $\Lambda $-packing problem is $NP$-hard, 
i.e. possibly intractable in general, this problem turns out 
to be tractable for some natural classes of graphs
(e.g. \cite{KKN,Kclfree}). 
It would be  also interesting to find polynomial algorithms 
that would provide a good approximation solution for 
the problem. Below (see {\bf \ref{km}}, {\bf \ref{2,3-graphs}},
and {\bf \ref{cubic-connected}})
are some examples of such results.
In each case the corresponding packing problem is 
polynomially solvable.

Let $v(G)$ and $\lambda (G)$ denote the number 
of vertices and the maximum number of disjoint
3--vertex paths in $G$, respectively.
Obviously $\lambda (G) \le \lfloor v(G)/3 \rfloor $.

In \cite{K,KM} we answered the following natural question:
\\[1ex]
\indent
{\em How many disjoint 3-vertex paths must a cubic $n$-vertex graph have?}
\bs 
\label{km} 
If $G$ is a cubic graph, then 
$\lambda (G) \ge \lceil v(G)/4  \rceil $.
Moreover, there is a polynomial time algorithm for 
finding a $\Lambda $-packing having at least  
$\lceil v(G)/4  \rceil$ components.
\es

Obviously if every component of $G$ is $K_4$, then
$\lambda (G) = v(G)/4$. 
Therefore the bound in {\bf \ref{km}} is sharp.
\\[.5ex]
\indent
Let ${\cal G}^3_2 $ denote the set of graphs with each vertex of degree at least $2$ and at most $3$.

In \cite{K} we  answered (among other results) the following  question:
\\[1ex]
\indent
{\em How many disjoint 3-vertex paths must an $n$-vertex  graph from ${\cal G}^3_2$ have?}

\bs 
\label{2,3-graphs} Suppose that $G \in {\cal G}^3_2$ and  
$G$ has no 5-vertex components.
Then $\lambda (G) \ge \lceil v(G)/4  \rceil $.
\es

Obviously {\bf \ref{km}}  follows from {\bf \ref{2,3-graphs}} because if $G$ is a cubic graph, then $G \in {\cal G}^3_2$ and $G$ has no 5-vertex components. 
\\[1ex]
\indent
In \cite{K} we  also gave a construction  that allowed  to prove the following:
\bs 
\label{extrgraphs1}
There are infinitely many connected graphs for  which the bound  in {\bf \ref{2,3-graphs}} is attained.
Moreover, there are infinitely many subdivisions of
cubic 3-connected graphs for which the bound  in 
{\bf \ref{2,3-graphs}} is attained.
\es

The next interesting question is:
\\[1ex]
\indent
{\em How many disjoint 3-vertex paths must a cubic connected graph have?}
\\[1ex]
\indent
In \cite{K2} we proved the following.
\bs 
\label{cubic-connected} 
Let ${\cal C}_n$ denote the set of connected cubic graphs with $n$ vertices and 
\\ 
$\lambda _n = \min \{\lambda (G)/v(G): G \in {\cal C}_n\}$.
Then for some $c > 0$,
%
%
\[\frac{3}{11}(1 - \frac{c}{n}) \le \lambda _n \le 
\frac{3}{11}(1 - \frac{1}{n^2}).\]
\es

The next natural question is:

\bs {\bf Problem}
\label{2con-must}
How many disjoint 3-vertex paths must a cubic 2-connected graph have?
\es

This question is still open  (namely, the sharp lower bound on the number  of disjoint 3-vertex paths in a cubic 2-connected $n$-vertex graph is unknown).

There are infinitely many 2-connected and cubic graphs 
such that
$\lambda (G) <   \lfloor v(G)/3 \rfloor $.
\\[1ex]
\indent
As to cubic 3-connected graphs, an old open questions here is:

\bs {\bf Problem}
\label{Pr3con} Is the following claim true:
\\[.5ex]
if $G$ is a 3-connected and cubic graph,
then $\lambda (G) =   \lfloor v(G)/3 \rfloor $ ?
\es

In \cite{K3con-cub} we discuss Problem {\bf \ref{Pr3con}} and
show, in particular, that the claim in 
{\bf \ref{Pr3con}} is equivalent to some seemingly 
much stronger claims. Here are some results of this kind.

\bs {\em \cite{K3con-cub}}
\label{cubic3-con}
The following are equivalent for cubic 3-connected graphs $G$:
\\[1ex]
${\bf (z1)}$ 
$v(G) = 0 \bmod 6$ $\Rightarrow$ $G$ has 
a $\Lambda $-factor,
\\[1ex] 
${\bf (z2)}$
$v(G) = 0 \bmod 6$ $\Rightarrow$ for every 
$e \in E(G)$ there is a $\Lambda $-factor of $G$ 
avoiding $e$,
\\[1ex]
${\bf (z3)}$
$v(G) = 0 \bmod 6$ $\Rightarrow$ for every 
$e \in E(G)$ there is a $\Lambda $-factor of $G$ 
containing $e$,
\\[1ex]
${\bf (z4)}$
$v(G) = 0 \bmod 6$ $\Rightarrow$ $G - X$ has 
a $\Lambda $-factor for every $X \subseteq E(G)$, $|X| = 2$,
\\[1ex]  
${\bf (z5)}$
$v(G) = 0 \bmod 6$ $\Rightarrow$ 
$G - L$ has a $\Lambda $-factor for every
3-vertex path $L$ in $G$,
\\[1ex]
${\bf (t2)}$
$v(G) = 2 \bmod 6$ $\Rightarrow$ $G - \{x,y\}$ 
has a $\Lambda $-factor for every $xy \in E(G)$,
\\[1ex]
${\bf (f1)}$
$v(G) = 4 \bmod 6$ $\Rightarrow$ $G - x$ 
has a $\Lambda $-factor for every $x \in V(G)$,
\\[1ex]
${\bf (f2)}$
$v(G) = 4 \bmod 6$ $\Rightarrow$ $G - \{x, e\}$ 
has a $\Lambda $-factor for every $x \in V(G)$ and 
$e \in E(G)$.
\es

In \cite{Kclfree} we have shown (in particular) that all claims in  {\bf \ref{cubic3-con}} except for $(z5)$ are true for 3-connected claw-free graphs and $(z5)$ is true for 
cubic, 2-connected, and claw-free graphs distinct from $K_4$ (see also  \cite{KKN}).

The problems similar to {\bf \ref{Pr3con}} are interesting for  2-connected and cubic 
graphs having some additional properties. For example,

\bs {\bf Problem} 
\label{Pr2con,cub,bip,pl}  
Is $\lambda (G) =   \lfloor v(G)/3 \rfloor $ true for every 
2-connected, cubic, bipartite, and planar  graph ?
\es

In this paper (see Section \ref{2con}) we answer the question in 
{\bf \ref{Pr2con,cub,bip,pl}} by 
giving a construction that provides infinitely many   
2-connected, cubic, bipartite, and planar  graphs such that  
$\lambda (G) <   \lfloor v(G)/3 \rfloor $
(see also \cite{K2con-cbp}).

\section{Some notation, constructions, 
 and simple observations}
\label{constructions}

\indent

\begin{figure}
  \centering
  \includegraphics{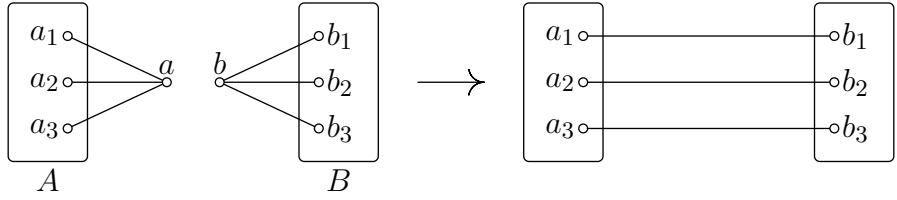}
  \caption{$Aa \sigma bB$}
  \label{fAasbB}
\end{figure}

We consider undirected graphs with no loops and 
no parallel edges unless stated otherwise. 
As usual, $V(G)$ and $E(G)$ denote the set of vertices and edges of $G$, respectively, and $v(G) = |V(G)$.
If $X$ is a vertex subset or a subgraph of $G$, then let
$D(X,G)$ or simply $D(X)$, denotes the set of edges in $G$, having exactly one end--vertex in $X$, and let
$d(X,G) = |D(X,G)|$.
If $x \in V(G)$, then $D(x,G)$ is the set of edges in $G$ incident to $x$,
$d(x,G) = |D(x,G)|$, $N(x,G) = N(x)$ is the set of vertices 
in $G$ adjacent to $x$, and 
$\Delta (G) = \max \{d(x,G): x \in V(G)\}$.
If $e = xy \in E(G)$, then let $End(e) = \{x,y\}$.
Let $Cmp (G)$ denote the set of components of $G$ and $cmp(G) = |Cmp(G)|$.
\\[1ex]
\indent
Let ${\cal C}(k)$ denote the set of cubic $k$-connected graphs, $k \in \{1,2,3\}$.
\\[1ex]
\indent
Let $A$ and $B$ be disjoint graphs, 
$a \in V(A)$,  $b \in V(B)$, and
$\sigma : N(a,A) \to N(b,B)$ be a bijection.
Let $Aa \sigma bB$ denote the graph
$(A - a) \cup (B - b) \cup \{x\sigma (x): x \in N(a,A)\}$.
We usually assume that 
$N(a,A) = \{a_1, a_2, a_3\}$, $N(b,B) = \{b_1, b_2, b_3\}$, and  $\sigma (a_i) = b_i$ for $i \in \{1,2,3\}$
%
(see Fig. \ref{fAasbB}). 

We also say that {\em $Aa \sigma bB$ is obtained from $B$ 
by replacing vertex $b$ by $(A - a)$ according to $\sigma $}.

Let $B$ be a cubic graph and $X \subseteq V(B)$. 
Let $A(v)$, where $v \in X$, be a graph, 
$a^v$ be  a vertex of degree three in $A(v)$, and 
$A^v = A(v) - a^v$.
By using the above operation, we can build a graph
$G = B\{(A(v), a^v): v \in X\}$ by replacing each vertex $v$
of $B$ in $X$ by $A^v $ assuming that all $A(v)$'s are disjoint. Let $D^v = D(A^v,G)$.
For each $u \in V(B) \setminus X$ let $A(u)$ be the graph having exactly two vertices $u$, $a^u$ and exactly three parallel edges connecting $u$ and $a^u$.
Then $G = B\{(A(v), a^v): v \in X\} = 
B\{(A(v), a^v): v \in V(B)\}$.
If, in particular, $X = V(G)$ and each $A(v)$ is a copy of $K_4$, then $G$ is obtained from $B$ by replacing each vertex by a triangle.

Let $E' = E(G) \setminus \cup \{E(A^v):  v \in V(B)\}$.
Obviously, there is a unique bijection
$\alpha : E(B) \to E'$ such that if $uv \in E(B)$, then 
$\alpha (uv)$ is an edge in $G$ having one end-vertex in
$A^u$ and the other in $A^v$.

Let $P$ be a $\Lambda $-packing in $G$.
For $uv \in E(B)$, $u \ne v$, we write  $u \neg ^p v$ or simply,
$u \neg  v$,
if $P$ has a 3-vertex path $L$ such that $\alpha (uv) \in E(L)$ and $|V(A^u) \cap V(L)| = 1$.
Let $P^v$ be the union of components of $P$ that meet
$D^v$ in $G$.
\\[1ex]
\indent
Obviously
\bs
\label{AasbB} 
Let $k$ be an integer and $k \le 3$.
If $A$ and $B$ above are $k$-connected, cubic, 
bipartite, and planar graphs, then $Aa\sigma bB$ 
is also a $k$-connected, cubic, bipartite, and 
planar graph, respectively.
\es

From {\bf \ref{AasbB}} we have: 
\bs
\label{A(B)} 
Let $k$ be an integer and $k \le 3$.
If $B$  and each $A_v$ is a $k$-connected, cubic, 
bipartite, and planar graphs, then
$B\{(A_v, a_v): v \in V(B)\}$ is also a $k$-connected, 
cubic, bipartite, and planar graph, respectively.
\es

\begin{figure}
  \centering
  \includegraphics{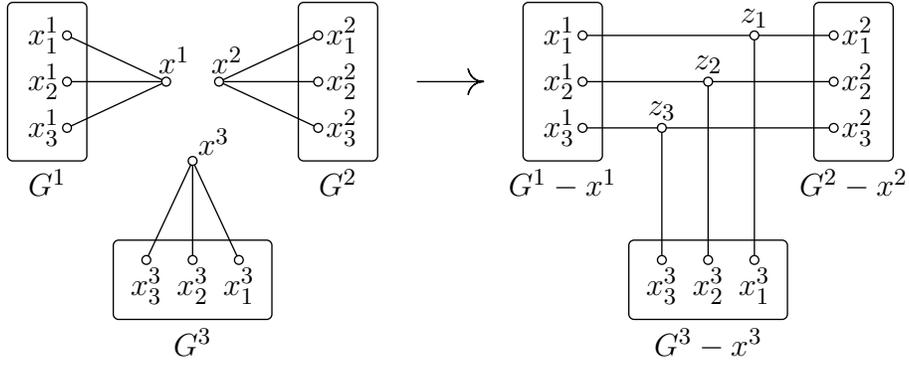}
  \caption{$Y(A^1,a^1; A^2,a^2; A^3,a^3)$}
  \label{fY}
\end{figure}

Let $A^1$, $A^2$, $A^3$ be three disjoint graphs,
$a^i \in V(A^i)$,  and $N(a^i,A^i) = \{a^i_1, a^i_2, a^i_3\}$, 
where $i \in \{1,2,3\}$.
Let  $F = Y(A^1,a^1; A^2,a^2; A^3,a^3)$  denote the graph 
obtained from $(A^1 - a^1) \cup (A^2-a^2) \cup (A^3 - a^3)$ 
by adding three new vertices $z_1$, $z_2$, $z_3$ 
and the set of nine new edges $\{z_ja^i_j:  i, j \in \{1,2,3\}$
(see Fig. \ref{fY}). 
In other words, if $B = K_{3,3}$ is the complete
$(X,Z)$-bipartite graph with 
$X = \{x_1,x_2,x_3\}$ and  $Z = \{z_1,z_2,z_3\}$, then $F$ is obtained from the $B$ by replacing each vertex $x_i$ in $X$ by $A^i - a^i$ so that
$D(A^i - a^i, F) = \{a^i_jz_j: j \in \{1,2,3\}$.
Let $D^i = D(A^i- a^i, F)$.
If $P$ is a $\Lambda $-packing of $F$, then let 
$P^i = P^i(F)$ be the union of components of $P$ meeting $D^i$ and $E^i = A^i(P) = E(P) \cap D^i$, $i \in \{1,2,3\}$.

If each $(A^i, a^i)$ is a copy of the same $(A, a)$, 
then we write $Y(A, a)$ instead of 
$Y(A^1, a^1; A^2,a^2; A^3,a^3)$.
\\

From {\bf \ref{A(B)} } we have, in particular: 
\bs
\label{Y(G1,G2,G3)} 
Let $k$ be an integer and $k \le 3$.
If each $A^i$ above is a $k$-connected, cubic, and 
bipartite graph, then $Y(A^1,a^1; A^2,a^2; A^3,a^3)$ 
$($see Fig. \ref{fY}$)$ is also a $k$-connected, cubic, and bipartite graph, 
respectively.
\es

Let $A$ and $B$ be disjoint graphs, 
$a = a_1a_2 \in E(A)$, and $b = b_1b_2 \in E(B)$.
Let $AabB$ be the graph obtained from 
$(A - a) \cup (B - b)$ by adding two new edges 
$a_1b_1$, $a_2b_2$
%
(see Fig. \ref{fAabB}).
\begin{figure}
  \centering
  \includegraphics{lambda-packing-5.mps}
  \caption{$AabB$}
  \label{fAabB}
\end{figure}
Let $Aa|bB$ be the graph obtained from $AabB$ 
by replacing edge $a_ib_i$ by a 3-vertex path 
$a_iz_ib_i$ for each $i \in \{1,2\}$ and by adding 
a new edge $z = z_1z_2$.
We call $z$ the {\em middle edge of} $Aa|bB$ 
%
(see Fig. \ref{fAa|bB}).
\begin{figure}
  \centering
  \includegraphics{lambda-packing-6.mps}
  \caption{$Aa|bB$}
  \label{fAa|bB}
\end{figure}
\\

It is easy to see the following.

\bs
\label{A|B} 
Let $k$ be an integer and $k \le 2$.
If $A$ and $B$ are $k$-connected, cubic, bipartite, 
and planar graphs, then
both $AabB$ and $Aa|bB$ are also a $k$-connected, 
cubic, bipartite, and planar graph, respectively.
\es

We will use the following simple observation.

\bs
\label{3cut}
Let $A$ and $B$ be disjoint graphs, 
$a \in V(A)$,  $N(a,A) = \{a_1, a_2, a_3\}$,  
$b \in V(B)$,  $N(b,B) = \{b_1, b_2, b_3\}$, and 
$G = Aa \sigma bB$, where each $\sigma (a_i) = b_i$
$($see Fig. \ref{fAasbB}$)$. 
Let $P$ be a $\Lambda $-factor of $G$ 
$($and so $v(G) = 0 \bmod 3$$)$ and 
$P'$ be the $\Lambda $-packing of $G$ consisting of the components $($3-vertex paths$)$ of $P$ that meet 
$\{a_1b_1, a_2b_2, a_3b_3\}$.
\\[1ex]
$(a1)$
Suppose that  $v(A) = 0 \bmod 3$, and 
so $v(B) = 2 \bmod 3$.  Then one of the following holds
$($see Fig \ref{A1}$)$:
\begin{figure}
  \centering
  \includegraphics{lambda-packing-3.mps}
  \caption{}
  \label{A1}
\end{figure}
\\[0.5ex]
\indent
$(a1.1)$ $P'$ has exactly one component that
has two  vertices in $A - a$, that are adjacent 
$($and, accordingly, exactly one vertex in $B - b$$)$,
\\[0.5ex]
\indent
$(a1.2)$ $P'$ has exactly two components  and each component has exactly one vertex in $A - a$ 
$($and, accordingly, exactly two vertices in 
$B - b$, that are adjacent$)$,
\\[0.5ex]
\indent
$(a1.3)$ $P'$ has exactly three components $L_1$, $L_2$, $L_3$ and one of them, say $L_1$, has exactly one vertex in $A - a$ and each of the other two $L_2$, $L_3$, has exactly two vertices in $A - a$, that are adjacent
$($and, accordingly, $L_1$ has exactly two  vertices in $B - b$, that are adjacent, and each of the other two $L_2$, $L_3$, has exactly  one vertex in $B - b$$)$. 
\\[1ex]
$(a2)$
Suppose that  $v(A) = 1 \bmod 3$, and 
so $v(B) = 1 \bmod 3$.  Then one of the following holds
$($see Fig \ref{A2}$)$:
\begin{figure}
  \centering
  \includegraphics{lambda-packing-4.mps}
  \caption{}
  \label{A2}
\end{figure}
\\[0.5ex]
\indent
$(a2.1)$ $P' = \emptyset $,
\\[0.5ex]
\indent
$(a2.2)$ $P'$ has exactly two components, say $L_1$, $L_2$, and one of the them, say $L_1$, has exactly one vertex in 
$A - a$ and exactly two  vertices in $B - b$, that are adjacent, and the other component $L_2$ has exactly two vertices in 
$A - a$, that are adjacent, and exactly one vertex in $B - b$,
\\[0.5ex]
\indent
$(a2.3)$ $P'$ has exactly three components $L_1$, $L_2$, $L_3$ and either each $L_i$ has exactly one vertex in 
$A - a$ 
$($and, accordingly, has exactly two vertices in 
$B - b$, that are adjacent$)$ or each $L_i$ has exactly two vertices in $A - a$, that are adjacent 
$($and, accordingly, has exactly one vertex in $B - b$$)$.
\es


\section{$\Lambda $-packings in cubic  2-connected  graphs}
\label{2con}

\indent

\subsection{Cubic, bipartite, and 2-connected graphs} 
\label{2con,cub,bip}

\bs
\label{Y,cmp(Pi)<3} 
Let $G = Y(A^1,a^1;  A^2,a^2; A^3,a^3)$ 
$($see Fig. \ref{fY}$)$ and 
$P$  be a $\Lambda $-factor of $G$.
Suppose that each $A^i$ is a cubic graph and 
$v(A^i) = 0 \bmod 6$.
Then $cmp (P^i) \in \{1, 2\}$ for every $i \in \{1,2,3\}$.
\es

\bp Let $i \in \{1,2,3\}$.
Since $D^i$ is a matching and $P^i$ consists of the components of $P$ meeting $D^i$, clearly $cmp (P^i) \le 3$.
Since $v(A^i) = -1 \bmod 6$, we have $cmp(P^i) \ge 1$.
It remains to show that $cmp (P^i) \le 2$.
Suppose, on the contrary, that $cmp(P^1) = 3$.

Since $P$ is a $\Lambda $-factor of $G$ and 
 $v(A^1 - a^1) = -1\bmod 6$, clearly
 $v(P^1) \cap V(A^1 - a^1) = 5$ and we can assume 
 (because of symmetry) that 
$P_1$ consists of three components  $a^1_3z_3a^2_3$,  
$z_1a^1_1 y^1$, and $z_2a^1_2 u^1$ for some 
$y^1, u^1 \in V(A^1)$. 
Then $cmp(P^3) = 0 $, a contradiction.
\ep

\bs
\label{Y,e-} 
Let $A$ be a graph, $e  = aa_1\in E(A)$, and 
$G = Y(A,a)$ $($see Fig. \ref{fY}$)$.
Suppose that
\\[0.5ex]
$(h1)$ $A$ is  cubic,
\\[0.5ex]
$(h2)$ $v(A) = 0\bmod 6$, and 
\\[0.5ex]
$(h3)$ $a$ has no $\Lambda $-factor containing 
$e = aa_1$.
\\[0.5ex]
\indent
Then $v(G) = 0 \bmod 6$ and $G$ has 
no $\Lambda $-factor.
\es

\bp (uses {\bf \ref{Y,cmp(Pi)<3}}).
Suppose, on the contrary, that $G$ has  a $\Lambda $-factor $P$. By definition of $G = Y(A,a)$, each $A^i$ is a copy of $A$ and edge $e^i = a^ia^i_1$ in $A^i$ is a copy of edge 
$e = aa_1$ in $A$. 
By {\bf \ref{Y,cmp(Pi)<3}}, $cmp(P^i) \in \{1, 2\}$
for every $i \in \{1,2,3\}$.
Since $P$ is a $\Lambda $-factor of $G$ and 
$v(A^i - x^i) = -1\bmod 6$, clearly
$E(P) \cap D^i $ is an edge subset of a $\Lambda $-factor 
of $A^i$ for every $i \in \{1,2,3\}$
(we assume that edge $z_ja^i_j$ in $G$ is edge $a^ia^i_j$ in $A^i$).
Since $a^1a^i_1$ belongs to no  $\Lambda $-factor of 
$A^i$ for every $i \in \{1,2,3\}$, clearly
$E(P) \cap \{z_1a^1_1,z_1a^2_1,z_1a^3_1\}  = \emptyset $.
Therefore $z_1 \not \in V(P)$, and so $P$ is not 
a $\Lambda $-factor of $G$, a contradiction.
\ep

\bs
\label{2conG,e-}
Suppose $A$ and $B$ are cubic graphs
and $v(A) = 2\bmod 6$,  $v(B) = 2\bmod 6$. 
Let $G = Aa|bB$ with the middle edge $z = (z_1z_2)$
$($see Fig. \ref{fAa|bB}$)$.
Then $v(G) = 0\bmod 6$ and $G$ has 
no $\Lambda $-factor containing edge $z$.
\es

\bp Obviously $v(G) = 0\bmod 6$.
Let $P$ be a $\Lambda $-factor of $G$.
Let $Z = \{z_1a_1,z_2a_2\}$. 
Since  $P$ is a $\Lambda $-factor of $G$ and 
$v(A) = 2\bmod 6$,
clearly $E(P) \cap R \ne \emptyset $.

Suppose that $|E(P) \cap Z| = 1$, say 
$E(P) \cap Z = z_1a_1$.
Since $P$ is a $\Lambda $-factor of $G$ and 
$v(A) = 2\bmod 6$, clearly $z_1a_1a \in P$ for some 
$a \in V(A)$. Then $z \not \in E(P)$.

Now suppose that $|E(P) \cap Z| = 2$.
Again since  $P$ is a $\Lambda $-factor of $G$ and 
$v(A) = 2\bmod 6$, clearly $a_iz_iy_i \in P$ for every 
$i\in \{1,2\}$ and some $y_1, y_2 \in V(A)$, $y_1 \ne y_2$. 
Since $y_1 \ne y_2$ and $z = z_1z_2 \in E(G)$, clearly 
$z \not \in E(P)$.
\ep
\\

A minimum simple cubic graph $H$ with $v(H) = 2 \bmod 6$
has 8 vertices. 
Therefore by {\bf \ref{2conG,e-}}, a minimum simple 
cubic graph $M$ with an edge $z$, avoidable by every $\Lambda $-factor of $M$, that can be obtained by construction $Aa|bB$, has 18 vertices

\begin{figure}
  \centering
  \includegraphics{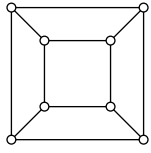}
  \caption{The cube}
  \label{fQ}
\end{figure}

\begin{figure}
  \centering
  \includegraphics{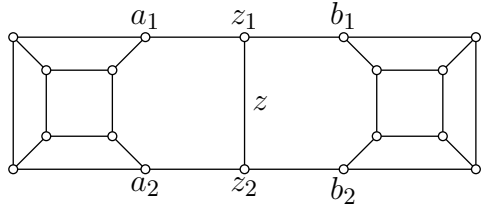}
  \caption{$K = Q_1q_1|q_2Q_2$, $v(K) = 18$}
  \label{fK}
\end{figure}

The graph-skeleton $Q$ of the cube is the only  simple, 
cubic, and bipartite  graph with 8 vertices 
%
(see Fig. \ref{fQ}). 
Moreover, $Q$ is planar. Therefore if both $A$ and $B$ are disjoint copies $Q_1$ and $Q_2$ 
of $Q$, then by  {\bf \ref{2conG,e-}}, 
$K = Q_1q_1|q_2Q_2$ has no $\Lambda $-factor, 
containing the middle edge $z$ of $K$, and  
$v(K) = 18$ 
%
(see Fig. \ref{fK}).
\begin{figure}
  \centering
  \includegraphics{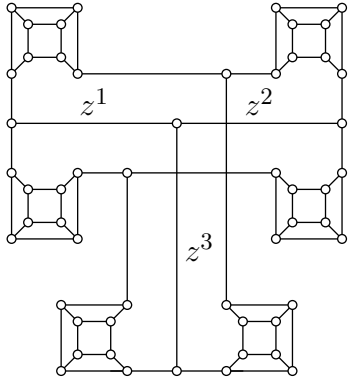}
  \caption{$R = Y(K,k)$, $v(R) = 54$}
  \label{fR}
\end{figure}
Moreover, since $Q$ is  2-connected, bipartite, and 
planar, $K$ is also 2-connected, bipartite, and planar.

\bs
\label{D(M,e)nofactor}
Suppose that $H$ is a cubic graph, $v(H) = 0 \bmod 6$, 
$h \in E(H)$, and $H$ has no $\Lambda $-factor containing $h$. Let $G = Y(H,h)$. 
Then $v(G) = 0 \bmod 6$ and $G$ has no  $\Lambda $-factor.
\es

\bp Obviously $(H,h)$ satisfies $(h1)$, $(h2)$, and $(h3)$ of 
{\bf \ref{Y,e-}}. 
Therefore by {\bf \ref{Y,e-}}, $G$ has no $\Lambda $-factor.
\ep 

\bs
\label{2con,nofactor}
There are infinitely many graphs $G$ such that $G$ is 
2-connected, cubic, and bipartite, $v(G) = 0 \bmod 6$, 
and $G$ has no $\Lambda $-factor.  
\es

\bp Follows immediately from {\bf \ref{AasbB}}, {\bf \ref{A|B}}, 
{\bf \ref{Y,e-}}, {\bf \ref{2conG,e-}}, and 
{\bf \ref{D(M,e)nofactor}}.
\ep
\\[1ex]
\indent
By {\bf \ref{D(M,e)nofactor}}, 
$R = Y(K,k))$ has no $\Lambda $-factor, where $k$ is a vertex in $K$ incident to the middle edge $z$ of $K$
(see Fig. \ref{fR}).
Obviously, $R$ is a 2-connected, cubic, bipartite graph 
and $v(R) = 54$. 
The graph $R$ is a smallest simple graph with these 
properties provided by the above construction.

\subsection{Cubic, bipartite, planar, and 2-connected graphs}
\begin{figure}
  \centering
  \includegraphics{lambda-packing-10.mps}
  \caption{}
  \label{fS}
\end{figure}
\begin{figure}
  \centering
  \includegraphics[width=3.4in]{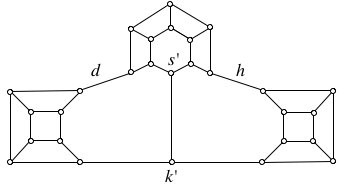}
  \caption{$H = Kk\sigma  sS$, $v(H) = 28$}
  \label{fH}
\end{figure}
\bs
\label{H,h+r-}
Let $A$ and $B$ be cubic graphs, 
$v(A) = 0 \bmod 6$ and $v(B) = 0 \bmod 6$,
$a \in V(A)$,  $N(a,A) = \{a_1, a_2, a_3\}$, 
$b \in V(B)$,  $N(b,B) = \{b_1, b_2, b_3\}$.
Let $H = Aa\sigma  bB$,
where $\sigma (a_i) = b_i$ $($see Fig. \ref{fAasbB}$)$.
Suppose that $A$ has no $\Lambda $-factor containing $aa_1$. 
Then $v(H) = 4 \bmod 6$ and
$H - b_2$ has no $\Lambda $-factor avoiding $h = a_3b_3$.
\es

\bp Let $P$ be a $\Lambda $-factor of $H- b_2$.
Since $v(A - a) = -1 \bmod 6$
and $A$ has no $\Lambda $--factor containing $aa_1$, 
clearly $b_3a_3z \in P$ for some vertex $z$ in $A$ 
adjacent to $a_3$, and so $h = a_3b_3 \in E(P)$.
\ep

\bs
\label{G,x-}
Suppose that $A$ and $B$ are disjoint cubic graph,
$v(A) = 0 \bmod 6$, $v(B) = 4 \bmod 6$,  
$a = a_1a_2 \in E(A)$, $x \in V(B)$, 
$b  = b_1b_2 \in E(B)$,  
$A$ has no $\Lambda $-factor containing $a$, and
$B - x$ has no $\Lambda $-factor avoiding $b$ 
$($and so $x$ is not incident to $b$$)$.
Let $G = AabB$ $($Fig \ref{fAabB}$)$. 
Then $v(G) = 4 \bmod 6$ and $G - x$ has 
no $\Lambda $-factor. 
\es

\bp Obviously $v(G) = 4 \bmod 6$.
Suppose, on the contrary, that $G - x$ has 
a $\Lambda $-factor $P$. Since $B - x$ has 
no $\Lambda $-factor and $v(G) = 0 \bmod 6$, clearly 
$\{a_1b_1, a_2b_2\} \subseteq E(P)$ and there are 
vertices $a_3 \in V(A)$ and $b_3 \in V(B)$ such that
(up to symmetry) $a_3a_1b_1, a_2b_2b_3 \in P$.
Then $(P \cap A) \cup a_3a_1a_2$ is a $\Lambda $-factor 
of $A$ containing $a$, a contradiction.
\ep

\bs
\label{G,e+}
Suppose that $A$ and $B$ are disjoint cubic graph,
$v(A) = 2 \bmod 6$, $v(B) = 4 \bmod 6$,  
$a = a_1a_2 \in E(A)$,  $b  = b_1b_2 \in E(B)$,  and
$B - b_1$ has no $\Lambda $-factor. 
Let $G = AabB$ $($Fig \ref{fAabB}$)$. 
Then $v(G) = 0 \bmod 6$ and $G$ has 
no $\Lambda $-factor avoiding $a_2b_2$. 
\es

\bp Suppose, on the contrary, that  $G$ has 
a $\Lambda $-factor $P$ avoiding $a_2b_2$.
Since $P$ is a $\Lambda $-factor of $G$ and
$v(A) = 2 \bmod 6$, clearly 
$E(P) \cap \{a_1b_1, a_2b_2\} = a_1b_1$ and 
$a_3a_1b_1 \in P$ for some $a_3 \in V(A - \{a_1,a_2\}$.
Then $P \cap B$ is a $\Lambda $-factor of $B - b_1$ 
a contradiction.
\ep

\bs
\label{Gnofactor}
Suppose that $A$ and $B$ are disjoint cubic graph,
$v(A) = 0 \bmod 6$, $v(B) = 0 \bmod 6$,  
$a = a_1a_2 \in E(A)$,  $b  = b_1b_2 \in E(B)$,  and
$A$ has no $\Lambda $-factor containing $a$, and
$B$ has no $\Lambda $-factor avoiding $b$. 
Let $G = AabB$ $($Fig \ref{fAabB}$)$. 
 Then $v(G) = 0 \bmod 6$ and $G$ has no 
$\Lambda $-factor. 
\es

\bp Suppose, on the contrary, that $G$ has 
a $\Lambda $-factor $P$. 
Since $P$ is a $\Lambda $-factor of $G$ and
$B$ has no $\Lambda $-factor avoiding $b$,
 clearly $ \{a_1b_1, a_2b_2\} \subseteq E(P) \cap$ 
and there are vertices 
$a_3 \in V(A)$ and $b_3 \in V(B)$ such that
(up to symmetry) $a_3a_1b_1, a_2b_2b_3 \in P$.
Then $(P \cap A) \cup a_3a_1a_2$ is a $\Lambda $-factor 
of $A$ containing $a$, a contradiction.
\ep

\bs
\label{2con,cub,bip,pl,nofactor}
There are infinitely many graphs $G$ such that
$v(G) = 0 \bmod 6$, $G$ is 2-connected, cubic, bipartite, 
and planar, and $G$ has no $\Lambda $-factor.
\es

\begin{figure}
  \centering
  \includegraphics[width=5.7in]{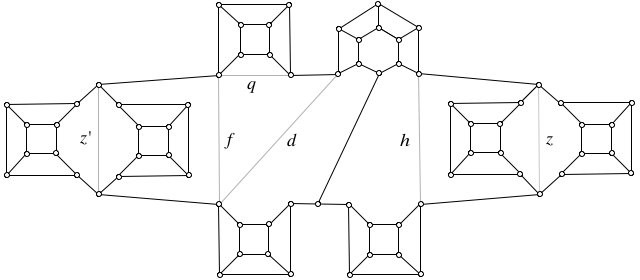}
  \caption{$N = Kz'fF$, $v(N) = 72$}
  \label{fN}
\end{figure}

\bp Follows from {\bf \ref{AasbB}}, {\bf \ref{A|B}}, 
{\bf \ref{2conG,e-}},  {\bf \ref{H,h+r-}}, {\bf \ref{G,x-}}, 
{\bf \ref{G,e+}}, and {\bf \ref{Gnofactor}}.
\ep
\\[2ex]
\indent
Let ${\cal CBP}$ denote the set of cubic, bipartite, planar, 
and 2-connected graphs.
The smallest  graph $G$ in ${\cal CBP}$ with 
$v(G) = 0 \bmod 6$ is the six-prism $S$, $v(S) = 12$
(see Fig. \ref{fS}).
\begin{figure}
  \centering
  \includegraphics[width=3.8in]{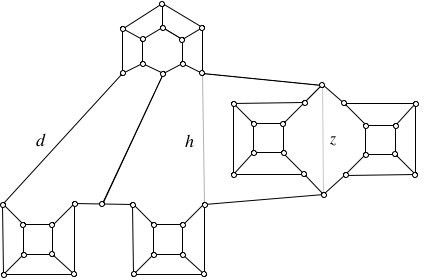}
  \caption{$D = KzhH$, $v(D) = 46$}
  \label{fD}
\end{figure}
Therefore by {\bf \ref{H,h+r-}}, the smallest  graph $H$ 
in ${\cal CBP}$ that has properties, guaranteed by  
{\bf \ref{H,h+r-}}, and that can be obtained by construction 
$Aa\sigma bB$ (see Fig \ref{fAasbB}), is $H = Kk\sigma  sS$, $v(H) = 28$ (see Fig \ref{fH}).
\begin{figure}
  \centering
  \includegraphics[width=3.88in]{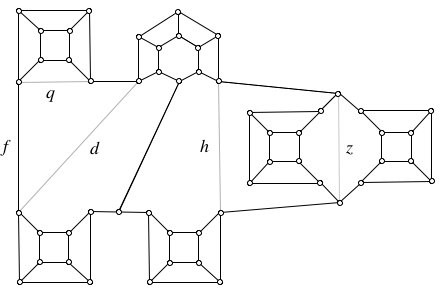}
  \caption{$F = Qqd(KzhH)$, $v(F) = 54$}
  \label{fF}
\end{figure}
Accordingly, by {\bf \ref{G,x-}}, the smallest  graph $D$ 
in ${\cal CBP}$ that has properties, guaranteed by  
{\bf \ref{G,x-}}, and that can be obtained by construction 
$AabB$ (see Fig \ref{fAabB}), is $D = KzhH$, $v(D) = 46$
(see Fig \ref{fD}).
%
%
By {\bf \ref{G,e+}}, the smallest  graph $F$ in
${\cal CBP}$ that has properties, guaranteed by  
{\bf \ref{G,e+}}, and can be obtained by construction 
$AabB$, is $F = QqdD$, $v(F) = 54$
(see Fig. \ref{fF}). 
Now by {\bf \ref{Gnofactor}}, the smallest  graph $N$ 
in ${\cal CBP}$ with $v(N) = 0  \bmod 6$ that has no 
$\Lambda $-factor, and that can be obtained by the above construction, is $N = Kz'fF$, $v(N) = 72$
%
(see Fig. \ref{fN}).

\end{document}